\documentclass[11pt]{amsart}

\theoremstyle{plain}
\newtheorem{thm}{Theorem}[section]
\newtheorem{theorem}[thm]{Theorem}

\newtheorem{lemma}[thm]{Lemma}
\newtheorem{corollary}[thm]{Corollary}
\newtheorem{proposition}[thm]{Proposition}
\theoremstyle{definition}
\newtheorem{remark}[thm]{Remark}

\newtheorem{definition}[thm]{Definition}

\newtheorem{example}[thm]{Example}

\newtheorem{conjecture}[thm]{Conjecture}

\numberwithin{equation}{section}

\newcommand{\bv}{{\mathbf v}}

\newcommand{\bX}{{\mathbf X}}
\newcommand{\bY}{{\mathbf Y}}
\newcommand{\bZ}{{\mathbf Z}}
\newcommand{\bH}{{\mathbf H}}

\newcommand{\p}{\partial}

\newcommand{\sD}{{\mathcal D}}
\newcommand{\sE}{{\mathcal E}}
\newcommand{\sF}{{\mathcal F}}

\newcommand{\sK}{{\mathcal K}}
\newcommand{\sL}{{\mathcal L}}

\newcommand{\sO}{{\mathcal O}}
\newcommand{\sP}{{\mathcal P}}

\newcommand{\sS}{{\mathcal S}}

\newcommand{\sU}{{\mathcal U}}

\newcommand{\sX}{{\mathcal X}}

\newcommand{\C}{{\mathbb C}}

\newcommand{\BP}{{\mathbb P}}

\newcommand{\fg}{{\mathfrak g}}

\newcommand{\fgl}{{\mathfrak g}{\mathfrak l}}

\newcommand{\fm}{{\mathfrak m}}

\newcommand{\fh}{{\mathfrak h}}

\title[Formal principle with convergence]{Formal principle with convergence for  rational curves of Goursat type}

\author{Jun-Muk Hwang}

\thanks{This work was supported by the Institute for Basic Science (IBS-R032-D1).}

\begin{document}
\begin{abstract}
We propose a conjecture that a general member of a bracket-generating family of rational curves in a complex manifold satisfies the formal principle with convergence, namely, any formal equivalence between such curves is convergent. If the normal bundles of the rational curves are positive, the conjecture follows from the results of Commichau-Grauert and Hirschowitz. We prove the conjecture for the opposite case when  the normal bundles are furthest from  positive vector bundles among bracket-generating families, namely, when the families of rational curves are of Goursat type.   The proof uses natural ODEs associated to rational curves of Goursat type and corresponding Cartan connections constructed by Doubrov-Komrakov-Morimoto.
As an example, we see that  a general line on a smooth cubic fourfold satisfies the formal principle with convergence.
\end{abstract}

\maketitle



\section{Introduction}\label{s.introduction}

In this paper, we use the following  definition of the formal principle.

\begin{definition}\label{d.FP} For a compact complex submanifold $A$ in a complex manifold $X$, denote by $(A/X)_{\infty}$ the formal neighborhood of $A$ in $X$ and by $(A/X)_{\sO}$ the germ of Euclidean neighborhoods of $A$ in $X$. We say that $A \subset X$
 satisfies {\em the formal principle},
    if for any  compact submanifold $\widetilde{A}$ in a complex manifold $\widetilde{X}$ admitting
      a formal isomorphism $(A/X)_{\infty} \stackrel{\cong}{\to} (\widetilde{A}/\widetilde{X})_{\infty}$ between the formal neighborhoods, there exists a biholomorphic isomorphism
$(A/X)_{\sO} \stackrel{\cong}{\to} (\widetilde{A}/\widetilde{X})_{\sO},$ namely,  a biholomorphic map between some Euclidean neighborhoods $U \subset X$ of $A$ and $\widetilde{U} \subset \widetilde{X}$ of $\widetilde{A}$ that sends $A$ to $\widetilde{A}$, giving the following commutative diagram.
$$ \begin{array}{ccccccc}   X & \supset & U & \stackrel{\cong}{\longrightarrow} & \widetilde{U} & \subset  & X \\ & & \cup & & \cup & & \\ & & A & \stackrel{\cong}{\longrightarrow} & \widetilde{A}& & \end{array}$$
 \end{definition}

We are interested in  the following stronger form of the formal principle.

\begin{definition}\label{d.convergence}
A compact complex submanifold $A$ in a complex manifold $X$ satisfies {\em the formal principle with convergence},  if for  any  compact submanifold $\widetilde{A}$ in a complex manifold $\widetilde{X}$ admitting
      a formal isomorphism $\psi: (A/X)_{\infty} \stackrel{\cong}{\to} (\widetilde{A}/\widetilde{X})_{\infty}$, there exists a biholomorphic isomorphism
$\Psi: (A/X)_{\sO} \stackrel{\cong}{\to} (\widetilde{A}/\widetilde{X})_{\sO}$
  such that $\psi= \Psi|_{(A/X)_{\infty}},$ giving the following commutative diagram.
  $$ \begin{array}{ccccccc}   X & \supset & U & \stackrel{\Psi}{\longrightarrow} & \widetilde{U} & \subset  & X \\ & & \cup & & \cup & & \\ & & (A/X)_{\infty} & \stackrel{\psi}{\longrightarrow} & (\widetilde{A}/\widetilde{X})_{\infty}. & & \end{array}$$
  In other words, any formal isomorphism $\psi$ is convergent.
\end{definition}

Definition \ref{d.convergence} is a much stronger property of $A \subset X$ than Definition \ref{d.FP}.
For example, when $A$ is a point in a complex manifold $X$, it satisfies the formal principle, but does not satisfy the formal principle with convergence.
 In this paper, we restrict our discussion to the case when $A$ is a smooth rational curve.
Many questions on the formal principle are already very interesting and difficult  for smooth rational curves.
 We use the following notion.

  \begin{definition}\label{d.positive} Let $C \subset X$ be a smooth rational curve. We say that a vector  bundle of rank $r$ on $C$  is {\em positive} if it is isomorphic to a direct sum of line bundles $\sO(a_1) \oplus \cdots \oplus \sO(a_{r})$ on $\BP^1 \cong C$ with $a_i >0$  for all $i$.
  We say that a vector bundle $V$ on $C$  is {\em semipositive} if it is isomorphic to a direct sum $V^+ \oplus \sO_C^{q}$ of a positive vector subbundle $V^+ \subset V$ and a complementary trivial bundle $\sO_C^q, q= {\rm rank} (V) - {\rm rank} (V^+),$ on $C$.
  \end{definition}

 Regarding the formal principle for smooth rational curves, the following result is proved in \cite{CG} and \cite{Hi}.

\begin{theorem}\label{t.CG} If $C \subset X$ is a smooth rational curve with positive normal bundle, then it satisfies the formal principle with convergence. \end{theorem}

On the other hand, the following result is proved in \cite{Hw19}.

\begin{theorem}\label{t.Hwang}
If $C \subset X$ is a smooth rational curve with semipositive normal bundle, then a general deformation of $C \subset X$ satisfies the formal principle. \end{theorem}

In Theorem \ref{t.Hwang}, a general deformation of $C$ does not necessarily satisfy  the formal principle with convergence. Here is an easy counterexample.

\begin{example}\label{e.product}
Let $C$ be a smooth rational curve in a complex manifold $Y$ with semi-positive normal bundle $N_{C/Y}.$
Let $X$ be the product $Y \times Z$ for a positive-dimensional complex manifold $Z$. Then $N_{C/X}$ is semipositive.  No deformation of $C$ in X can satisfy the formal principle with convergence, because there are many formal isomorphisms $(z/Z)_{\infty} \cong (z/Z)_{\infty}$ at a point $z \in Z$ which do not converge. \end{example}

We would like to formulate a suitable geometric condition to avoid examples like Example \ref{e.product} and strengthen Theorem \ref{t.Hwang}.  For this purpose, we introduce the following definitions.

\begin{definition}\label{d.distribution}
A {\em distribution} $D$ on a complex manifold $M$ is a saturated subsheaf $D \subset \Theta_M$ of the tangent sheaf of $M$. Denote by $\p^1 D \subset \Theta_M$ the distribution spanned by $[D, D] + D$ and inductively define $\p^{i+1} D$ as the distribution  spanned by $ [\p^i D, \p^i D] + \p^i D$ with $\p^0 D = D$.  A distribution $D$ is {\em bracket-generating} if $\p^i D = \Theta_M$ for some $i$. \end{definition}

\begin{definition}\label{d.sD}
For a complex manifold $X$, consider the Douady space ${\rm Douady}(X)$ of all compact analytic subspaces of $X$ (see Section VIII.1 of \cite{GPR} for an overview with references).
Let $\sS(X) \subset {\rm Douady}(X)$ be the space of smooth rational curves on $X$ with semipositive normal bundles. Let $\sK$ be a connected open subset in $\sS(X).$ Then $\sK$ is a complex manifold and for a smooth rational curve $C \subset X$ corresponding to a point $[C] \in \sK$, the tangent space $T_{[C]}\sK$ is naturally isomorphic to the vector space $H^0(C, N_{C/X})$ of holomorphic sections of the normal bundle $N_{C/X}$ of $C \subset X$. Denote by $\sD_{[C]}$ the subspace of $T_{[C]} \sK$ corresponding to $$H^0(C, N^+_{C/X}) \subset H^0(C, N_{C/X}),$$ which is independent of the choice of an isomorphism $$N_{C/X} \cong N_{C/X}^+ \oplus \sO_C^q, \ q= {\rm rank} (N_{C/X}) - {\rm rank} (N^+_{C/X})$$ for the semipositive vector bundle $N_{C/X}$ in the notation of Definition \ref{d.positive}. This determines a natural distribution $\sD \subset \Theta_{\sK}.$ We say that $\sK$ is a {\em bracket-generating family of rational curves} if $\sD$ is a bracket-generating distribution on $\sK$. \end{definition}

 It is clear that $\sD = \Theta_{\sK}$ if and only if the normal bundles of members of $\sK$ are positive. Thus we can restate Theorem \ref{t.CG} as follows.

\begin{theorem}\label{t.CG2}
If $\sD =\Theta_{\sK}$ in Definition \ref{d.sD}, then all members of $\sK$ satisfy the formal principle with convergence. \end{theorem}

It is easy to see that deformations of $C$ in Example \ref{e.product} do not belong to a bracket-generating family. We propose the following conjecture, as a generalization of Theorem \ref{t.CG2} and a strengthening of Theorem \ref{t.Hwang}.

\begin{conjecture}\label{c.FPC}
Let $\sK$ be a bracket-generating family of rational curves on a complex manifold. Then a general member of $\sK$ satisfies the formal principle with convergence. \end{conjecture}

In other words, if the positive deformation directions generate all directions by successive Lie brackets, then we expect the formal principle with convergence for general members.
In this paper, we prove a special case of Conjecture \ref{c.FPC}. To state our result, we recall the following.

\begin{definition}\label{d.Goursat}
 A distribution $D \subset \Theta_M$ on a complex manifold $M$  is a {\em Goursat} distribution if the rank of $\p^i D$ is $i+2$ for all $0 \leq i \leq \dim M -2$.
A bracket-generating family $\sK$ of rational curves on a complex manifold $X$ is a {\em family of rational curves of Goursat type} if the distribution $\sD \subset \Theta_{\sK}$ is a Goursat distribution.\end{definition}

\begin{remark}\label{r.Goursat}
In Definition \ref{d.Goursat}, a Goursat distribution has rank 2. There is a more general definition of  Goursat distributions allowing higher ranks (Definition 6.4 of \cite{Mon}). The study of such high rank cases can often be reduced to the case of rank 2. \end{remark}

Our main result is the following (see Theorem \ref{t.main2} for a more precise formulation).

\begin{theorem}\label{t.main}
Let $\sK$ be a family of rational curves of Goursat type on a complex manifold.
Then a general member of $\sK$ satisfies the formal principle with convergence. \end{theorem}

When $\sD$ is a Goursat distribution, we have $\dim H^0(C, N^+_{C/X}) = 2$ for any $[C] \in \sK$. This implies $$N_{C/X} \cong \sO(1) \oplus \sO^{\dim X -2}.$$
Thus $N_{C/X}$ is the least positive among nontrivial semipositive bundles. Moreover, the Goursat condition ${\rm rank}(\p^i \sD) = i+2$ for all $0 \leq i \leq \dim \sK-2$ means that the growth of Lie brackets of $\sD$ is the slowest possible among bracket-generating distributions. For these reasons, one may say that the normal bundles of rational curves of Goursat type have the weakest positivity property among families of bracket-generating rational curves. Thus Theorem \ref{t.main} shows that Conjecture \ref{c.FPC} holds even when the normal bundles of the rational curves have the weakest positivity property. We believe that this is a good evidence toward Conjecture \ref{c.FPC}.

Are there many examples of rational curves of Goursat type? When $\dim M \leq 4$, any bracket-generating distribution of rank 2 on $M$ is a Goursat distribution. This provides many interesting examples of rational curves of Goursat type on complex manifolds of dimension $4$. For example, general lines on smooth cubic hypersurface in $\BP^5$ are of Goursat type (Example 4.9 of \cite{HL}), which gives the following corollary of Theorem \ref{t.main}.

\begin{corollary}\label{c.cubic} A general line on a smooth cubic hypersurface in $\BP^5$ satisfies the formal principle with convergence. \end{corollary}

See Example \ref{e.HL} below for many examples of rational curves of Goursat type in higher dimensions.

The proof of Theorem \ref{t.main} uses the idea of Cartan's equivalence method as in \cite{Hw19}. This approach to the study of the formal principle is based on the fact that  formal equivalences of rational curves induce formal equivalences of
the natural differential systems at the corresponding points in the deformation spaces of the curves. Using this approach, we may say that the proof of Theorem \ref{t.Hwang} in \cite{Hw19} has the following two components, the geometric part to obtain a suitable involutive system ({\bf I.S.}) and the analytic part of applying Cartan-K\"ahler theorem ({\bf C.K.}):
\begin{itemize} \item[({\bf I.S.})] the natural differential system on the deformation space of rational curves leads to a suitable involutive system   by \cite{Mo83}; \item[({\bf C.K.})] the analytic solvability of  formally solvable involutive systems  follows from  Cartan-K\"ahler theorem. \end{itemize}
 For the proof of Theorem \ref{t.main}, we replace ({\bf I.S.}) by  a more refined geometric result of  Doubrov-Komrakov-Morimoto (\cite{DKM}), a construction of canonical Cartan connections or absolute parallelisms. Once we have natural Cartan connections, the analytic part is just the convergence of formal equivalences of affine connections proved in \cite{KN}. In other words, the proof of Theorem \ref{t.main} has the following two components:
 \begin{itemize} \item[({\bf A.P.})] the natural differential system on the deformation space of the rational curves leads to a canonical absolute parallelism by  \cite{DKM}; \item[({\bf K.N.})] the convergence of formal equivalences of absolute parallelisms follows from  Kobayashi-Nomizu theorem (Theorem \ref{t.KN} below). \end{itemize}

It is expected that this method of combining ({\bf A.P.}) with ({\bf K.N.}) can be used to check Conjecture \ref{c.FPC} in many cases other than rational curves of Goursat type.  From this perspective,  it seems valuable to construct canonical absolute parallelisms associated to natural  geometric structures on the deformation spaces, for various  bracket-generating families of rational curves.

Ngaiming Mok has asked us whether there is an example where some member of $\sK$ in Theorem \ref{t.main} does not satisfy the formal principle with convergence. We expect that such an example exists, although we have not yet been able to find one.

\section{Jet spaces and ODE-structures }
In this section, we introduce the notion of ODE-structures. This is a special case of a G-structure subordinate to a distribution studied in \cite{Mo93}.

\begin{definition}\label{d.prD}
Let $D \subset \Theta_M$ be a distribution on a complex manifold $M$ given by  a vector subbundle of the tangent bundle $TM$, which we denote by the same symbol $D \subset TM$ by abuse of notation.
\begin{itemize}
\item[(i)] The {\em Cauchy characteristic} of $D$ is the subdistribution ${\rm Ch}(D) \subset D$ spanned by local sections $v$ of $D$ satisfying $[v, D] \subset D.$
\item[(ii)]
        For a holomorphic submersion $\eta: \widetilde{M} \to M$ from a complex manifold $\widetilde{M}$, let ${\rm d} \eta: T\widetilde{M} \to TM$ be its derivative and let
        $\eta^{-1}D \subset T\widetilde{M}$ be the vector subbundle whose fiber at a point $\alpha \in \widetilde{M}$ is $({\rm d}\eta)^{-1} D_y, y = \eta(\alpha)$. The distribution on $\widetilde{M}$ given by $\eta^{-1}D$ is the {\em inverse image} of $D$ under $\eta$, which we denote by the same symbol $\eta^{-1} D \subset \Theta_{\widetilde{M}}$. It is immediate that ${\rm Ker}({\rm d} \eta) \subset {\rm Ch}(\eta^{-1}D).$
        \item[(iii)]
Let  $\eta: \BP D \to M$ be the projectivization of the vector bundle $D$.  For each $\alpha \in \BP D$ corresponding to a 1-dimensional subspace $\widehat{\alpha} \subset D_y$ with $y=\eta(\alpha),$
        denote by ${\rm pr} D_{\alpha} \subset (\eta^{-1}D)_{\alpha}$ the subspace $({\rm d} \eta)^{-1}(\widehat{\alpha}).$ This determines a vector subbundle ${\rm pr}D$ of $\eta^{-1}D$ with ${\rm rank} ({\rm pr} D) = {\rm rank} (D)$. The corresponding distribution ${\rm pr} D \subset \Theta_{\BP D}$ is the {\em prolongation } of the distribution $D$. \end{itemize}
\end{definition}

We recall the following from Section 2.3 of \cite{HL}.

\begin{definition}\label{d.jet}
Let $\C_t$ (resp. $\C_u$) be the complex line with the independent variable $t$ (resp. the dependent variable $u$) as an affine coordinate. Let $J^0 = \C_t \times \C_u$ be the $2$-dimensional complex manifold with the natural projection $p^0: J^0 \to \C_t$. Let $W^1 = TJ^0$ be the tangent bundle and let $E^1 \subset W^1$ be the line subbundle given by ${\rm Ker} ({\rm d} p^0)$. Define a sequence of $\BP^1$-bundles inductively
$$J^0 \stackrel{p^1}{\longleftarrow} \BP W^1 \stackrel{p^2}{\longleftarrow} \BP W^2 \longleftarrow \cdots\longleftarrow \BP W^{k-1} \stackrel{p^k}{\longleftarrow} \BP W^k \longleftarrow \cdots, $$
where
  $W^{k+1}$ is the prolongation of $W^k$, which is a vector subbundle of rank 2 in $T \BP W^k$ and contains the line subbundle $E^{k+1} \subset T \BP W^k$ defined  by ${\rm Ker} ({\rm d} p^{k})$.
    Define $J^1 := \BP W^1 \setminus \BP E^1$ and inductively define $$J^k := (\BP W^k \cap (p^k)^{-1}(J^{k-1})) \setminus \BP E^k.$$  We call $J^k$ the {\em  space of $k$-jets of functions} from $\C_t$ to $\C_u$. Writing $p^k|_{J^k}$ simply as $p^k$, we obtain a sequence of  affine $\C$-bundles
    $$J^0 \stackrel{p^1}{\longleftarrow} J^1 \stackrel{p^2}{\longleftarrow} J^2 \longleftarrow \cdots \longleftarrow  J^{k-1} \stackrel{p^k}{\longleftarrow}  J^k \longleftarrow \cdots. $$
    Regarding ${\rm d}t$ and ${\rm d} u$ as holomorphic functions on $TJ^0 = W^1$ which are homogeneous along the fibers of $W^1 \to J^0$, we have the meromorphic function $$u^{(1)} := \frac{{\rm d}u}{{\rm d}t}$$ on $\BP W^1$ which is a holomorphic function on $J^1$. Inductively, we have the meromorphic function $$u^{(k+1)} := \frac{{\rm d}u^{(k)}}{{\rm d} t}$$ on $\BP W^{k+1}$, which is a holomorphic function on $J^{k+1}$.
    It is clear that $(t, u, u^{(1)}, \ldots, u^{(k)})$ give a holomorphic coordinate system on the $(k+2)$-dimensional complex manifold $J^k$.
\end{definition}

We skip the proof of the following  easy lemma.

\begin{lemma}\label{l.jet} For each $k \geq 1,$
the distribution $W^{k+1}$ on $\BP W^k$ is a Goursat distribution satisfying $$\p^1 W^{k+1} = (p^k)^{-1} W^k \mbox{ and }  {\rm Ch}(\p^1 W^{k+1}) = E^{k+1}.$$ \end{lemma}

It is convenient to introduce the following terminology.

\begin{definition}\label{d.ODE}
Let $\sU$ be a complex manifold of dimension $k+2 \geq 4.$ An ordered pair of two line  subbundles $\sL, \sF \subset T \sU$ is an {\em ODE-structure} on $\sU$,  if  each point $u \in \sU$ has a neighborhood $u \in U \subset \sU$ equipped with  a biholomorphic map $\iota: U \stackrel{\cong}{\to} \iota(U) $ to an open subset $\iota(U) \subset J^k$ satisfying $${\rm d} \iota (\sF) = E^{k+1}|_{\iota(U)} \mbox{ and } W^{k+1}|_{\iota(U)} = {\rm d} \iota(\sF|_U) \oplus {\rm d} \iota(\sL|_U).$$ \end{definition}

The following elementary fact (e.g. Section 3.1 of \cite{DKM}) explains the motivation for the name `ODE-structure'. This lemma is necessary only for Example \ref{e.HL}.

    \begin{lemma}\label{l.tfae}
    For each $k \geq 2$, there are  natural one-to-one correspondences  among the following three classes of objects. \begin{itemize}
    \item[(a)]
    A germ of holomorphic sections of the holomorphic submersion $p^{k+1}: J^{k+1} \to J^k$ in a neighborhood of a point in $J^k$.
    \item[(b)] A germ of  ODEs of order $k+1$ of the form $$u^{(k+1)} = F(t, u, u^{(1)}, \ldots, u^{(k)})$$ in the independent variable $t$ and the dependent variable $u$, where $F$ is a holomorphic function  in a neighborhood of a point in $J^k$. \item[(c)] A germ of holomorphic foliations of rank 1 in a neighborhood $U$ of a point in $J^k$ given by a line subbundle  $\sL \subset W^{k+1}|_U \subset TU$ which gives a splitting   $$   W^{k+1}|_U = E^{k+1}|_U \oplus \sL.$$ \end{itemize} \end{lemma}

\section{Canonical Cartan connection associated to an ODE-structure}

In this section, we review the result on  the canonical Cartan connection associate to an ODE-structure. When $k=2$, an ODE-structure is a special case of parabolic geometries and the canonical Cartan connection is explicitly given in \cite{SY} extending the earlier work of S.-S. Chern \cite{Ch}. It follows also from  Tanaka's general construction in  \cite{Ta}.  When $k \geq 3$, an ODE-structure is no longer a parabolic geometry, and the canonical Cartan connection is constructed  in \cite{DKM} (see also  Section 2 of \cite{Do08} for a  succinct summary).  Let us mention that the presentation of \cite{DKM} and \cite{Do08} are in the real differentiable setting, but the translation to our holomorphic setting is immediate. We recall the definition of a Cartan connection first.

\begin{definition}\label{d.connection}
Let $G$ be a complex Lie group and let $H \subset G$ be a closed subgroup. Let $\fh \subset \fg$ be their Lie algebras.
A {\em Cartan connection} modelled after the homogeneous space $G/H$ on a complex manifold $M$ is an $H$-principal bundle $\pi: \sP \to M$ equipped with a $\fg$-valued 1-form $\omega$ on $\sP$ such that \begin{itemize}
\item[(i)] $\omega(\vec{v}) = v $ for the fundamental vector field $\vec{v}$ on $\sP$ associated to each $v \in \fh$;
    \item[(ii)] the right action $R_h: \sP \to \sP$ by any $h\in H$ satisfies $R_h^* \omega = {\rm Ad}(h)^{-1} \omega$; and
        \item[(iii)] $\omega$ defines an absolute parallelism on $\sP$, namely, an isomorphism $\omega_{\alpha}: T_{\alpha} \sP \cong \fg$ at each point $\alpha \in \sP$. \end{itemize} \end{definition}

We are interested in Cartan connections modelled after the following  homogeneous space $G/H$, the complex version of the homogeneous space in Section 2.2 of \cite{Do08}.

\begin{definition}\label{d.fg}
For $k \geq 2$, let $V_{k+1}$ be the vector space of homogeneous polynomials of degree $k+1$ in two variables $x, y$. We choose the following basis of $V_{k+1}$.
$$\{ \bv_i := \frac{x^{k+1-i} y^i}{i !} \mid 0 \leq i \leq k+1 \}.$$
Let $\fgl_2$ be the Lie algebra spanned by four vector fields $$\bX := x \frac{\p}{\p y}, \ \bY := y \frac{\p}{\p x}, \ \bH:= x \frac{\p}{\p x} - y \frac{\p}{\p y}, \ \bZ:= x \frac{\p}{\p x} + y \frac{\p}{\p y}.$$ Let ${\rm GL}_2$ be the simply connected complex Lie group with Lie algebra $\fgl_2$. The natural $\fgl_2$-representation on $V_{k+1}$ gives an action of ${\rm GL}_2$ on $V_{k+1}$. Let $G$ be the semidirect product of ${\rm GL}_2$ and $V_{k+1}$. Its Lie algebra $\fg$ has the following gradation. \begin{eqnarray*} \fg_1 &=& \C \bY \\ \fg_0 & =& \C \bH + \C \bZ \\ \fg_{-1} & = & \C \bX + \C \bv_{k+1} \\ \fg_{-i} & = & \C \bv_{k+2-i} \mbox{ for } 2 \leq i \leq k+2. \end{eqnarray*}
Let $\fh$ be $\fg_1 + \fg_0$ and let $H \subset G$ be the closed subgroup with Lie algebra $\fh \subset \fg$.
\end{definition}

The following is essentially Corollary in Section 6.4 of \cite{DKM} (also Proposition 2 of \cite{Do08}).

\begin{theorem}\label{t.DKM}
Let $G/H$ be as in Definition \ref{d.fg}.
Let $\sU$ be a complex manifold of dimension $k+2 \geq 4$ with an ODE-structure  $\sL, \sF \subset T \sU$.
Then there exists a Cartan connection $(\pi: \sP \to \sU, \omega)$ modelled after $G/H$  canonically associated with the ODE-structure.
\end{theorem}

Note that in \cite{DKM} (resp. \cite{Do08}), an ODE-structure  is called a $G_0$-structure of type $\fm$ (resp.  a geometric structure inherited from the jet space): the line bundle $E^k$ in Section 3.3 of \cite{DKM} (resp.  the line bundle $E$ in Section 2.1 of \cite{Do08}) is our $\sL$ and the vertical distribution $V^k$ in Section 3.3 of \cite{DKM} (resp. the line bundle $F$ in Section 2.1 of \cite{Do08}) is our $\sF$.

As mentioned before, the proof of Theorem \ref{t.DKM} is an immediate translation of the proof of Corollary in Section 6.4 of \cite{DKM} into holomorphic setting. Thus instead of giving a full proof, we just give some comments. The proof in \cite{DKM} used Morimoto's theory of filtered manifolds, which could look rather technical for those without previous knowledge of the theory. We recommend Sections 6.1 and 6.2  of \cite{DKM} for a brief summary of Morimoto's theory.  The key point of the proof of Theorem \ref{t.DKM} uses Morimoto's criterion (Proposition 3.10.1 in \cite{Mo93} or Theorem 2 in \cite{DKM}) for the construction of canonical Cartan connections. This criterion is a condition on the algebraic property of the graded Lie algebra $\fg.$ In the original setting of \cite{DKM}, the corresponding graded Lie algebra  is a real graded Lie algebra  defined in the same way as in Definition \ref{d.fg} with real coefficients instead of complex coefficients. For this real graded Lie algebra,  Morimoto's criterion
 has been checked in the proof of Theorem 4 in \cite{DKM} (or Section 2.3 of \cite{Do08}) using the positive definite inner product on the Lie algebra with the orthogonal basis $\bX, \bY, \bH, \bZ, \bv_i,$ satisfying  $$ \langle \bX, \bX \rangle = \langle \bY, \bY \rangle = 1, \ \langle \bH, \bH \rangle = \langle \bZ, \bZ \rangle = 2, $$ $$ \langle \bv_i, \bv_i\rangle = (n-i)!/i! \  \mbox{ for } \ 0 \leq i \leq k+1.$$ In our case where $\fg$ is a complex Lie algebra, we can define $\langle \cdot, \cdot \rangle$ as a positive definite Hermitian inner product with the above properties. Then all arguments of \cite{DKM} work verbatim to give a proof of Theorem \ref{t.DKM} in the holomorphic setting.

\medskip
  By the functorial nature of the construction in \cite{Mo93}, the Cartan connection in Theorem \ref{t.DKM} is canonically associated with the ODE-structure in  the following sense. Let
$\widetilde{\sU}$ be another complex manifold with an ODE-structure  $\widetilde{\sL}, \widetilde{\sF} \subset T \widetilde{\sU}$ and  the associated Cartan connection $(\widetilde{\pi}: \widetilde{\sP} \to \widetilde{\sU}, \widetilde{\omega}).$
Suppose that  there exists a biholomorphic map $\Phi: U \stackrel{\cong}{\to} \widetilde{U}$ which preserves the ODE-structures,   satisfying $${\rm d} \Phi(\sL) = \widetilde{\sL} \mbox{ and } {\rm d} \Phi (\sF) = \widetilde{\sF}.$$
Then it induces a biholomorphic map
$ \Phi_*: \sP \stackrel{\cong}{\to} \widetilde{\sP}$ satisfying the commutative diagram
$$ \begin{array}{ccc} \sP & \stackrel{\Phi_*}{\longrightarrow} & \widetilde{\sP} \\
\pi \downarrow & & \downarrow \widetilde{\pi} \\ \sU & \stackrel{\Phi}{\longrightarrow} & \widetilde{\sU} \end{array} $$
and $ \omega = \widetilde{\omega} \circ {\rm d} \Phi_*.$

Let us briefly sketch this construction, risking  oversimplification.   We have a sufficiently large integer $k >0$ such that the principal bundle $\sP$ (resp. $\widetilde{\sP}$) can be described as a fiber subbundle of ${\rm Frame}^k (\sU)$ (resp. ${\rm Frame}^k (\widetilde{\sU})$),  the bundle  of $k$-th order frames (in a suitable sense) on $\sU$ (resp. $\widetilde{\sU}$).    Any biholomorphic map $\Phi: U \to \widetilde{U}$ induces a biholomorphic map $\Phi^k: {\rm Frame}^k (\sU)
\to {\rm Frame}^k (\widetilde{\sU})$. For a general biholomorphism $\Phi$, we cannot have $\Phi^k(\sP) = \widetilde{\sP}$, but if we require that $\Phi$ preserves the ODE-structures, then $\Phi^k(\sP) = \widetilde{\sP}$ holds, inducing the above commutative diagram. 

 Similarly, suppose that there exists a
 formal isomorphism $\varphi: (u/\sU)_{\infty} \stackrel{\cong}{\to} (\widetilde{u}/\widetilde{\sU})_{\infty}$  for some points $u \in \sU$ and $\widetilde{u} \in \widetilde{\sU}$ satisfying
$$ {\rm d} \varphi (\sL|_{(u/\sU)_{\infty}}) = \widetilde{\sL}|_{(\widetilde{u}/\widetilde{\sU})_{\infty}} \mbox{ and }  {\rm d} \varphi (\sF|_{(u/\sU)_{\infty}}) = \widetilde{\sF}|_{(\widetilde{u}/\widetilde{\sU})_{\infty}}. $$ Then it induces a formal isomorphism
 $$\varphi_*: (\pi^{-1}(u)/\sP)_{\infty} \stackrel{\cong}{\to} (\widetilde{\pi}^{-1}(\widetilde{u})/ \widetilde{\sP})_{\infty}$$ satisfying the commutative diagram
$$ \begin{array}{ccc} (\pi^{-1}(u)/\sP)_{\infty} & \stackrel{\varphi_*}{\longrightarrow} &  (\widetilde{\pi}^{-1}(\widetilde{u})/ \widetilde{\sP})_{\infty} \\ \pi \downarrow & & \downarrow \widetilde{\pi} \\
 (u/\sU)_{\infty} & \stackrel{\varphi}{\longrightarrow} &  (\widetilde{u}/\widetilde{\sU})_{\infty} \end{array}$$ and
 \begin{equation}\label{e.omega} \omega|_{(\pi^{-1}(u)/\sP)_{\infty}} = \widetilde{\omega} \circ {\rm d} \varphi_*.\end{equation}
 
 In fact,  for a positive integer $i > > k$, we can choose a biholomorphic map $\varphi_i: U_i \to \widetilde{U}_i$ between open neighborhoods $u \in U_i \subset \sU$ and $\widetilde{u} \in \widetilde{U}_i \subset \widetilde{\sU}$ such that $\varphi_i$ agrees with $\varphi_{\infty}$ up to order $i$ at $u$. Then the ODE-structure ${\rm d} \varphi_i (\sL), {\rm d} \varphi_i (\sF) \subset T\widetilde{U}_i$ agrees with the ODE-structure $\widetilde{\sL}|_{\widetilde{U}_i}, \widetilde{\sF}|_{\widetilde{U}_i} \subset T\widetilde{U}_i$ up to order $i-1$ along the fiber over $\widetilde{u}$. Under the induced biholomorphic map $\varphi_i^k: {\rm Frame}^k (U_i) \to {\rm Frame}^k(\widetilde{U}_i),$ the image $\varphi_i^k(\sP)$ agrees with $\widetilde{\sP}$ along the fiber $\pi^{-1}(u)$ up to order $i-1-k$.  The $(i-1-k)$-th order neighborhood of the fiber $\widetilde{\pi}^{-1}(\widetilde{u})$  in the image $\varphi_i^k(\sP)$ does not depend on the choice of the biholomorphic map $\varphi_i$.    By taking $i$ arbitrarily large, we obtain the induced formal isomorphism $\varphi_*$.

\section{ODE-structures arising from rational curves of Goursat type}\label{s.Goursat}

In this section, we show that a family of rational curves of Goursat type gives rise to a complex manifold with an ODE-structure in a natural way.  First, we have the following result, which holds for families of rational curves somewhat more general than those of Goursat type.

\begin{proposition}\label{p.Univ}
Let $X$ be a complex manifold of dimension $k+1$. Let $\sK$ be a connected open subset in $\sS(X)$ in Definition \ref{d.sD} such that each point $[C] \in \sK$ corresponds to a smooth rational curve $C \subset X$ whose normal bundle $N_{C/X}$ is isomorphic to $\sO(1) \oplus \sO^{\oplus (k-1)}. $  In particular, the dimension of $\sK$ is $k+1$ and  the rank of the vector subbundle $\sD \subset T\sK$  in Definition \ref{d.sD} is 2.
Let $\rho: \BP \sD \to \sK$ be the projectivization of $\sD$. Then there exist a submersion $\mu: \BP \sD \to X$ such that for each $[C] \in \sK$, the morphism $\mu|_{\rho^{-1}([C])}:  \rho^{-1}([C]) \to X$ is an embedding whose image is the rational curve $C \subset X$ corresponding to $[C] \in \sK$. In other words, the double fibration
$$ \sK \stackrel{\rho}{\longleftarrow} \BP \sD \stackrel{\mu}{\longrightarrow} X$$ is isomorphic to the  universal family morphisms $$\sK \stackrel{\rho_{\sK}}{\longleftarrow} {\rm Univ}_{\sK} \stackrel{\mu_{\sK}}{\longrightarrow} X$$ associated with the open subset $\sK$ in the Douady space ${\rm Douady}(X)$. Moreover, the line subbundles ${\rm Ker}({\rm d} \mu)$ and ${\rm Ker}({\rm d} \rho)$ of $T \BP \sD$ satisfy
\begin{equation}\label{e.zero}
{\rm Ker}({\rm d} \mu) \cap {\rm Ker}({\rm d} \rho) =0 \mbox{ and  } \end{equation}
\begin{equation}\label{e.sum}  {\rm Ker}({\rm d} \mu) + {\rm Ker}({\rm d} \rho) = {\rm pr} \sD. \end{equation}
\end{proposition}

\begin{proof}
By standard deformation theory of rational curves  (see Section 3 of \cite{HL} and the references therein), for each $x \in X$, the subscheme $\sK_x $ in $\sK$ consisting of members of $\sK$ passing through $x$ is a submanifold of $\sK$ and its tangent space $T_{[C]}\sK_x$ at $[C] \in \sK_x $ can be identified with the vector subspace   $$  H^0(C, N_{C/X} \otimes {\bf m}_x) \subset H^0(C, N_{C/X}) = T_{[C]}\sK$$ consisting of holomorphic sections of $N_{C/X}$ vanishing at $x$. From \begin{eqnarray}\label{e.mx}
T_{[C]}\sK_x = H^0(C, N_{C/X} \otimes {\bf m}_x) & \subset & H^0(C, N^+_{C/X}) = \sD_{[C]}, \end{eqnarray}
the one-dimensional vector space $T_{[C]}\sK_x$ is contained in $\sD_{[C]}$ and corresponds to a point in  $\BP \sD_{[C]}$, which we denote by $\chi(x, C) \in \BP \sD_{[C]}$.
This determines a morphism $\chi: {\rm Univ}_{\sK} \to \BP \sD$.
From (\ref{e.mx}) with $x$ varying on $C$, we see that
$\chi$ sends the fiber $\rho_{\sK}^{-1}([C])$ isomorphically to the fiber of $\rho^{-1}([C]) = \BP \sD_{[C]}$. It follows that $\chi$ is a biholomorphic map satisfying the commutative diagram $$\begin{array}{ccccc}
\sK & \stackrel{\rho}{\longleftarrow} & \BP \sD & \stackrel{\mu}{\longrightarrow} & X \\ \| & & \chi \downarrow \cong & & \| \\ \sK & \stackrel{\rho_{\sK}}{\longleftarrow} & {\rm Univ}_{\sK} & \stackrel{\mu_{\sK}}{\longrightarrow} & X.\end{array} $$
Since $\mu$ embeds each fiber of $\rho$ as a smooth rational curve in $X$, we see (\ref{e.zero}). (\ref{e.sum}) follows from (\ref{e.mx}) and the definition of $\chi$.
\end{proof}

We recall the following classical result of Cartan (Theorem 6.5 of \cite{Mon}) on the structure of a Goursat distribution at a general point.

\begin{lemma}\label{l.Cartan}
  Let $D$ be a Goursat distribution on a complex manifold $M$ of dimension $k+1$.
  Write $\p^{(1)} D = \p^1 D$ and define the distribution $ \p^{(i+1)} D \subset \Theta_M$ inductively as the one spanned by $ [D, \p^{(i)} D] + \p^{(i)} D$ (where $\p^{(i)} D$ is written as $D^{i+1}$ in Section 2.3 of \cite{Mon}).
  Let $M' \subset M$ be the  Zariski-open subset consisting of points $y \in M$ satisfying the following two conditions.
  \begin{itemize}
    \item[(i)] For each $i \geq 0$, the distribution $\p^i D$  comes from  a vector subbundle of $TM$ in a neighborhood of $y$.  \item[(ii)]  $\p^i D = \p^{(i)} D$ for all $i \geq 0$  in a neighborhood of $y$. \end{itemize}
   Then any point $y \in M'$ admits a neighborhood $y \in O \subset  M'$ with a biholomorphic map $\zeta: O \to \zeta(O) \subset J^{k-1}$ satisfying ${\rm d} \zeta (D|_O) = W^{k}|_{\zeta(O)}.$ \end{lemma}

Proposition \ref{p.Univ} applied to a family of rational curves of Goursat type gives rise to an ODE-structure in the following way.

\begin{proposition}\label{p.ODE}
In Proposition \ref{p.Univ}, suppose that $\sK$ is a family of rational curves of Goursat type.   Let $\sK' \subset \sK$ be the Zariski-open subset obtained from Lemma \ref{l.Cartan} applied to the Goursat distribution $\sD \subset T \sK$. Then
for each point $y \in \sK'$, we have       a Euclidean neighborhood $y \in O \subset \sK'$ and a hypersurface $\sE \subset \BP \sD|_O$  such that on the complex manifold $\sU :=  \BP \sD|_{O} \setminus \sE,$  the line subbundles $\sL := {\rm Ker}({\rm d} \mu|_{\sU})$ and $ \sF := {\rm Ker}({\rm d} \rho|_{\sU})$ of $ T\sU$  determine an ODE-structure on $\sU$. \end{proposition}

\begin{proof}
By  the definition of $\sK'$,  any point $y \in \sK'$ admits a neighborhood $y \in O \subset \sK'$ with  a biholomorphic map $\zeta: O \to \zeta(O) \subset J^{k-1}$ satisfying ${\rm d} \zeta (\sD|_O) = W^{k}|_{\zeta(O)}.$ Let  $$ \eta: U:= \BP \sD|_O \stackrel{\cong}{\longrightarrow} \BP W^{k}|_{\zeta(O)} =: \eta(U)$$ be the biholomorphic map  induced by ${\rm d} \zeta.$  Then ${\rm d} \eta$ sends ${\rm pr} \sD|_{U }$ to $${\rm pr} W^k |_{\eta(U)} = W^{k+1} |_{\eta(U)} $$ and ${\rm Ker} (
{\rm d} \rho|_U)$ to $E^{k+1}$. By (\ref{e.zero}) and  (\ref{e.sum}), it sends the line subbundle ${\rm Ker}({\rm d} \mu|_U)$ of ${\rm pr}\sD|_U$  to a line subbundle of  $W^{k+1}|_U,$ which is transversal to $E^{k+1}|_U$. Put $\sE:= \eta^{-1} (\BP E^k|_{\zeta(O)})$ and $\sU := U \setminus \sE$ such that  $\eta(\sU) \subset J^k.$
Then  $\sL := {\rm Ker}({\rm d} \mu|_{\sU})$ and $\sF:={\rm Ker} ({\rm d} \rho|_{\sU})$
determine an ODE-structure on  $\sU$.   
 \end{proof}

We close this section with the following example of rational curves of Goursat type.

\begin{example}\label{e.HL}
By Theorem 1.1 in \cite{HL}, the ODE $$u^{(k+1)} = F(t, u, u^{(1)}, \ldots, u^{(k)})$$ associated to   the ODE-structure of Proposition \ref{p.ODE} via Lemma \ref{l.tfae} is of the form \begin{equation}\label{e.cubic}
u^{(k+1)} = a_3 (u^{(k)})^3 + a_2 (u^{(k)})^2 + a_1 u^{(k)} + a_0 \end{equation}
where $a_3, a_2, a_1, a_0$ are holomorphic functions in the $k+1$ variables $$(t, u, u^{(1)}, \ldots, u^{(k-1)})$$ defined in some domain in $\C^{k+1}$. Conversely, given any ODE of the type (\ref{e.cubic}), there exist
\begin{itemize} \item[(i)] a $(k+1)$-dimensional complex manifold $X$; \item[(ii)]  a family $\sK$ of rational curves of Goursat type  on $X$; and \item[(iii)] a biholomorphic map $\nu$ from a germ $\sX$ of the space of solutions of (\ref{e.cubic}) to an open subset $\nu(\sX)$ of $X$ \end{itemize} such that for any solution $u_0(t), t \in \Delta$ belonging to $\sX$, if we denote by $u^s_c(t), s \in \Delta,$ the solution satisfying  the initial conditions at $t=s$,  $$u(s) = u_0(s), u^{(1)}(s) = u_0^{(1)}(s), \ldots, u^{(k-1)}(s) = u_0^{(k-1)}(s), u^{(k)}(s) = c$$   for $c \in \C$ in a neighborhood $\Delta' \subset \C$ of $u_0^{(k)}(s)$, then the germ of the curve $$\{ [u^s_c] \in \sX \mid c \in \Delta'\}$$ is sent to the germ of a member of $\sK$ in $\nu(\sX) \subset X$. This shows that any ODE of the form (\ref{e.cubic}) yields an example  of a family of rational curves of Goursat type in a $(k+1)$-dimensional complex manifold.
 \end{example}

\section{Proof of Theorem \ref{t.main}}
We use the following result due to Kobayashi-Nomizu.

\begin{theorem}\label{t.KN}
Let $(M, \nabla)$ and $(\widetilde{M}, \widetilde{\nabla})$ be two complex manifolds with holomorphic affine connections.  Suppose there exist points $y \in M$ and $\widetilde{y} \in \widetilde{M}$ with a formal isomorphism $\xi: (y/M)_{\infty} \to (\widetilde{y}/\widetilde{M})_{\infty}$ such that $\xi$ sends the restriction  $\nabla|_{(y/M)_{\infty}}$ to the restriction $\widetilde{\nabla}|_{ (\widetilde{y}/\widetilde{M})_{\infty}}.$ Then $\xi$ is convergent, namely,  there exists a biholomorphic map $\Xi: (y/M)_{\sO} \to (\widetilde{y}/\widetilde{M})_{\sO}$ such that $\xi= \Xi|_{(y/M)_{\infty}}$ and $\Xi$ sends the restriction  $\nabla|_{(y/M)_{\sO}}$ to the restriction $\widetilde{\nabla}|_{ (\widetilde{y}/\widetilde{M})_{\sO}}.$ \end{theorem}

\begin{proof}
The assumption implies that $\xi$ sends all successive covariant derivatives of the torsion and the curvature of $\nabla$ at $y$ to those of $\widetilde{\nabla}$ at $\widetilde{y}$. Thus by Theorem 7.2 and its proof in Chapter VI of \cite{KN}, there is an isomorphism $\Xi$ of the affine connections from a neighborhood of $y$ to a neighborhood of $\widetilde{y}$ extending $\xi$. \end{proof}

We prove the following refined version of
Theorem \ref{t.main}.

\begin{theorem}\label{t.main2}
Let $\sK$ be a family of rational curves of Goursat type on a complex manifold
with the Goursat distribution $\sD$ on $\sK$. Let $\sK' \subset\sK$ be the closed analytic proper subset defined in Proposition \ref{p.ODE}. Then any member of $\sK \setminus \sK'$ satisfies the formal principle with convergence. \end{theorem}

\begin{proof}
Let $C \subset X$ be a member of $\sK \setminus \sK'$.
Suppose there is a formal isomorphism $\psi: (C/X)_{\infty} \to (\widetilde{C}/\widetilde{X})_{\infty}$ to a rational curve $\widetilde{C}$ in a complex manifold $\widetilde{X}$. Then $N_{C/X} \cong N_{\widetilde{C}/\widetilde{X}}$ and $[\widetilde{C}] \in  \sS(\widetilde{X})$  (in the notation of Definition \ref{d.sD}). Let $\widetilde{\sK} \subset  \sS(\widetilde{X})$ be a neighborhood of the point $[\widetilde{C}] \in \sS(\widetilde{X})$ with the natural distribution $\widetilde{\sD} \subset T \widetilde{\sK}$.

As explained in page 509 of \cite{Hi} (see also Lemma  3.5 of \cite{Hw19}), the functoriality of $X \mapsto {\rm Douady}(X)$  and the associated universal family implies that the formal isomorphism $\psi$
induces a formal isomorphism $$\psi_*: ([C]/\sK)_{\infty} \stackrel{\cong}{\longrightarrow}  ([\widetilde{C}]/\widetilde{\sK})_{\infty}$$ which sends $\sD|_{([C]/\sK)_{\infty}} $ to $\widetilde{\sD}|_{ ([\widetilde{C}]/\widetilde{\sK})_{\infty}}$  by the identification of $\BP \sD$ with the universal family ${\rm Univ}_{\sK}$ in Proposition \ref{p.Univ}.
Moreover, its derivative induces a formal isomorphism \begin{equation}\label{e.psistar} {\rm d} \psi_*: (\rho^{-1}([C])/ \BP \sD)_{\infty} \stackrel{\cong}{\longrightarrow} (\widetilde{\rho}^{-1}([\widetilde{C}])/\BP \widetilde{\sD})_{\infty}\end{equation} that sends ${\rm Ker}({\rm d} \mu)$ to ${\rm Ker}({\rm d} \widetilde{\mu})$ by the identification of $\mu$ (resp. $\widetilde{\mu}$) with the universal family morphism $\mu_{\sK}$ (resp. $\widetilde{\mu}_{\widetilde{\sK}}$) in Proposition \ref{p.Univ}.

 From the definition of $\sK'$ in Lemma \ref{l.Cartan}, the distribution  $\widetilde{\sD}$ is also a Goursat distribution and $\widetilde{\sK}$ is a family of rational curves of Goursat type. Moreover, if $\widetilde{\sK}' \subset \widetilde{\sK}$ is the closed analytic subset determined by the Goursat distribution $\widetilde{\sD}$ as in Proposition \ref{p.ODE}, then $[\widetilde{C}] $ belongs to $ \widetilde{\sK} \setminus \widetilde{\sK}'$ by Lemma \ref{l.Cartan}. Thus we have a neighborhood $O \subset \sK \setminus \sK'$ (resp. $\widetilde{O} \subset \widetilde{\sK} \setminus \widetilde{\sK}'$) of $[C]$ (resp. $\widetilde{C}$) and a complex manifold $\sU \subset \BP \sD|_O$ (resp. $\widetilde{\sU} \subset \BP \widetilde{\sD}|_{\widetilde{O}}$)  as in Proposition \ref{p.ODE} equipped with the natural ODE-structure $\sL, \sF \subset T \sU$ (resp. $\widetilde{\sL}, \widetilde{\sF} \subset T \widetilde{\sU})$.  Let  $(\pi: \sP \to \sU, \omega)$ (resp. $(\widetilde{\pi}: \widetilde{\sP} \to \widetilde{\sU}, \widetilde{\omega})$) be the canonical Cartan connection  from Theorem \ref{t.DKM} associated to this ODE-structure on $\sU$ (resp. $\widetilde{\sU}$).

Let $u \in \rho^{-1}([C]) \cap \sU$ be a point and let $\widetilde{u} \in \widetilde{\rho}^{-1}([\widetilde{C}]) \cap \widetilde{\sU}$ be its image under
${\rm d} \psi_*$ in (\ref{e.psistar}).
The formal isomorphism ${\rm d} \psi_*$ induces a formal isomorphism
$$\varphi: (u /\sU)_{\infty} \stackrel{\cong}{\longrightarrow}  (\widetilde{u}/\widetilde{\sU})_{\infty}$$ preserving the ODE-structures.
By the canonicity in   Theorem \ref{t.DKM},
we have the induced formal isomorphism
$$\varphi_*: (\pi^{-1}(u)/\sP)_{\infty} \to (\widetilde{\pi}^{-1}(\widetilde{u})/\widetilde{\sP})_{\infty}$$
satisfying (\ref{e.omega}).
The absolute parallelism $\omega$ (resp. $\widetilde{\omega}$ ) is an affine connection on $\sP$ (resp. $\widetilde{\sP}$). Thus (\ref{e.omega}) implies that $\varphi_*$ is convergent by Theorem \ref{t.KN}. It follows that $\varphi$ is convergent and there is a biholomorphic map
$\Phi: (u/\sU)_{\sO} \to (\widetilde{u}/\widetilde{\sU})_{\sO}$ such that $\varphi= \Phi|_{(u/\sU)_{\infty}}$.  Since $\varphi$ sends the line bundles $\sL, \sF$ to $\widetilde{\sL}, \widetilde{\sF}$, respectively, so does $\Phi$. Then $\Phi$ descends to a biholomorphic map $\Psi_*: ([C]/\sK)_{\sO} \to ([\widetilde{C}]/\widetilde{K})_{\sO}$ such that $\Psi_*|_{([C]/\sK)_{\infty}} = \psi_*$. Since ${\rm d} \Psi_*$ sends $\sD|_{([C]/\sK)_{\sO} }$ to $\widetilde{\sD}|_{([\widetilde{C}]/\widetilde{K})_{\sO}}$, it induces a biholomorphic map $$\Psi': (\rho^{-1}([C])/\BP \sD)_{\sO} \stackrel{\cong}{\longrightarrow} (\widetilde{\rho}^{-1}([\widetilde{C}])/ \BP \widetilde{\sD})_{\sO},$$ which is compatible with $\rho, \widetilde{\rho}, \mu, \widetilde{\mu}.$
Hence $\Psi'$ descends to a biholomorphic map $\Psi: (C/X)_{\sO} \stackrel{\cong}{\to} (\widetilde{C}/\widetilde{X})_{\sO}$ such that $\psi= \Psi|_{(C/X)_{\infty}}$. This shows that $C \subset X$ satisfies the formal principle with convergence.   \end{proof}

\bigskip
{\bf Acknowledgment}
I would like to thank the referee for valuable suggestions to improve the presentation of the paper.

\bigskip
Jun-Muk Hwang(jmhwang@ibs.re.kr)

\smallskip

Center for Complex Geometry,
Institute for Basic Science (IBS),
Daejeon 34126, Republic of Korea

\bigskip

\end{document}